\documentclass[12pt]{article}

\usepackage{amsmath,amssymb,amsthm}
\usepackage{amsfonts}
\usepackage[mathscr]{eucal}
\begin{document}
\newtheorem{thm}{Theorem}
\numberwithin{thm}{section}
\newtheorem{lemma}[thm]{Lemma}
\newtheorem{remark}{Remark}
\newtheorem{corr}[thm]{Corollary}
\newtheorem{proposition}{Proposition}
\newtheorem{theorem}{Theorem}[section]
\newtheorem{deff}[thm]{Definition}
\newtheorem{case}[thm]{Case}
\newtheorem{prop}[thm]{Proposition}
\numberwithin{equation}{section}
\numberwithin{remark}{section}
\numberwithin{proposition}{section}
\newtheorem{corollary}{Corollary}[section]
\newtheorem{others}{Theorem}
\newtheorem{conjecture}{Conjecture}\newtheorem{definition}{Definition}[section]
\newtheorem{cl}{Claim}
\newtheorem{cor}{Corollary}
\newcommand{\ds}{\displaystyle}
%\date{}

\newcommand{\stk}[2]{\stackrel{#1}{#2}}
\newcommand{\dwn}[1]{{\scriptstyle #1}\downarrow}
\newcommand{\upa}[1]{{\scriptstyle #1}\uparrow}
\newcommand{\nea}[1]{{\scriptstyle #1}\nearrow}
\newcommand{\sea}[1]{\searrow {\scriptstyle #1}}
\newcommand{\csti}[3]{(#1+1) (#2)^{1/ (#1+1)} (#1)^{- #1
 / (#1+1)} (#3)^{ #1 / (#1 +1)}}
\newcommand{\RR}[1]{\mathbb{#1}}
\thispagestyle{empty}
\begin{titlepage}
\title{\bf  Iterated Brownian Motion in Parabola-Shaped Domains}
\author{Erkan Nane\thanks{ Supported in part by NSF Grant \# 9700585-DMS }\\
Department of Mathematics\\
Purdue University\\
West Lafayette, IN 47906 \\
enane@math.purdue.edu}
\maketitle
\begin{abstract}
\noindent {\it Iterated Brownian motion $Z_{t}$ serves as a physical model
for diffusions in a crack. If $\tau _{D}(Z) $ is the first exit time of this processes   from a domain  $D \subset \RR{R}^{n}$, started at $z\in D$, then $P_{z}[\tau _{D}(Z)  >t]$ is the distribution
 of the lifetime of the process in $D$.  In this paper we
determine the large time asymptotics of $P_{z}[\tau _{P_{\alpha}}(Z)  > t]$ which gives  exponential
integrability of $\tau _{P_{\alpha}}(Z) $ for parabola-shaped domains of
the
form $ P_{\alpha}=\{ (x,Y)\in \RR{R}
\times \RR{R}^{n-1}: x>0,  |Y|<Ax^{\alpha} \}$, for $  0<\alpha <1$, $A>0.$ We also
obtain similar results for twisted domains in $\RR{R}^{2}$  as defined in \cite{DSmits}. In
particular, for a planar  iterated Brownian motion in a parabola $\mathcal{P}=\{(x,y):\ x>0,
 \ |y|< \sqrt{x} \}$ we find that for
$z\in \mathcal{P}$
$$ \lim_{t\to\infty}  t^{-\frac{1}{7}}  \log   P_{z}[\tau _{\mathcal{P}}(Z) >t]=
- \frac{7 \pi ^{2}}{2^{25/ 7}}. $$}
\end{abstract}

\end{titlepage}

\section{ Introduction}

Iterated Brownian motion has been of considerable interest
 to several authors in recent years; see for example
 Burdzy \cite{burdzy1, burdzy2},
DeBlassie \cite{deblassie}, Koshnevisan and Lewis
\cite{klewis} and references in these articles. Although this processes is not a
Markov process (it does not  satisfy the Chapman--Kolmogorov equations), it does have connections with the
 parabolic operator
$\frac{1}{8} \Delta ^{2} -\frac{\partial}{\partial t}$, as described in
  Funaki \cite{funaki}
 and DeBlassie \cite{deblassie}.

In analogy with ordinary Brownian motion and diffusions, if $\tau _{D}(Z)$ is the
first exit time of iterated Brownian motion from domain $D
$, started at $z\in D$,
 $P_{z}[\tau _{D}(Z) > t]$ provides a measure of the lifetime of the process in $D$.
The tail distribution of $\tau _{D}(Z)$ has a double integral representation in terms
of the probability density function of the  Brownian motion, as given in \cite{deblassie}.
This representation can then be used to
 compute the  asymptotics of the tail distribution of $\tau _{D}(Z)$ from which one can then obtain the
sharp  order of integrability  of  $\tau _{D}(Z)$. The goal of this paper is to do exactly this when the domain is
a general parabola in $\RR{R}^n$, or a twisted domain in $\RR{R}^{2}$.

To define the iterated  Brownian motion $Z_{t}$ started at $z \in \RR{R}$,
 let $X_{t}^{+}$, $X_{t}^{-}$ and $Y_{t}$
 be three  independent
one-dimensional
Brownian motions, all started at $0$. Two-sided Brownian motion is defined to be
\[ X_{t}=\left\{ \begin{array}{ll}
X_{t}^{+}, &t\geq 0\\
X_{(-t)}^{-}, &t<0.
\end{array}
\right. \]
Then the iterated Brownian motion started at $z \in \RR{R}$ is
\[ Z_{t}=z+X(Y_{t}),\ \ \    t\geq 0.\]
In $\RR{R}^{n}$, one requires $X^{\pm}$ to be independent
$n-$dimensional Brownian motions. This is the version of the iterated Brownian motion due to 
Burdzy. Our choice of independent
one-dimensional Brownian motions is motivated by the pde connection to
$\frac{1}{8} \Delta ^{2} -\frac{\partial}{\partial t}$ by Funaki
\cite{funaki} and DeBlassie \cite{deblassie} and by the interpretation
of the process as a diffusion in a crack by Burdzy and Khoshnevisan
\cite{bukh}.

The path properties of this process have been studied by Burdzy in  \cite{burdzy1}, and
\cite{burdzy2}. His works in particular imply the  LIL(Law of Iterated Logarithm) for the iterated Brownian motion.
 Khoshnevisan and Lewis extend results in \cite{burdzy2}, to
develop a stochastic calculus for iterated Brownian motion. The local time for
the iterated Brownian motion has been studied by Cs\`{a}ki, Cs\"{o}rg\H{o}, F\"{o}ldes and R\'{e}v\'{e}sz
\cite{csaki}, Eisenbum and Shi \cite{eisenbaumshi} and Xiao \cite{xiao}.

In \cite{deblassie}, DeBlassie obtains large time asymptotics of the tail
distribution of the exit time of $Z_t$  in bounded domains which have
regular boundary and for general cones in $\RR{R}^{n}$.
Let $D$ be a domain in $\RR{R}^{n}$.  Let
$$
\tau _{D} (Z) =\inf \{t\geq 0: \  Z_{t} \notin D \}
$$be the first exit time of $Z_{t}$ from $D$.
In this paper we prove the large time asymptotics of the tail distribution of $\tau _{D}(Z)$ when the
domain $D$  is a general parabola in $\RR{R}^{n}$, or a twisted domain in $\RR{R}^{2}$.  These type of questions
are now well understood in the case
of Brownian motion \cite{bds}, \cite{vandenberg}, \cite{DSmits}, \cite{lshi} and serve, together with DeBlassie's results in \cite{deblassie},  as motivation for our results.
Indeed,  Ba\~{n}uelos, DeBlassie and Smits showed in
\cite{bds} that if
$\tau _{\mathcal{P}}$ is the first exit time of the
Brownian  motion from the parabola $\mathcal{P}=\{(x,y):
x>0, |y|<A\sqrt{x} \}$, $A>0$, then there exist positive constants
$A_{1}$ and $A_{2}$ such that for $z\in
\mathcal{P}$
\begin{eqnarray}
-A_{1} \ \  & \leq  &\ \ \liminf_{t \rightarrow \infty} \  t^{-\frac{1}{3}}\ \log \ P_{z}[\tau _{\mathcal{P}} >t] \nonumber \\
& \leq &\ \ \limsup_{t \rightarrow \infty}\  t^{-\frac{1}{3}}\ \log \ P_{z}[\tau _{\mathcal{P}} >t] \leq \ \ -A_{2}. \nonumber
\end{eqnarray}
More recently, Lifshits and Shi \cite{lshi} found that the above limit exists for parabolic regions in any
dimension.  Let $0 <\alpha <1$ and $A>0$. We define the parabola-shaped domains as
$P_{\alpha}=\{ (x,Y)\in \RR{R} \times \RR{R}^{n-1}:
x>0, |Y|<Ax^{\alpha} \}$. Let $\tau _{\alpha}$ be the first exit time of Brownian motion from
$P_{\alpha}$.
Lifshits and Shi \cite{lshi}
proved that
\begin{equation}
\lim_{t\to\infty}t^{-(\frac{1-\alpha } { 1+\alpha} )}
\log  P_{z}[\tau _{\alpha}> t]=- l,
\end{equation}
where
\begin{equation}\label{brownianlimit}
l=(\frac{1+\alpha}{ \alpha})\left( \frac{\pi \jmath_{(n-3)/2}^{2/ \alpha}}{A^{2}2^{(3\alpha+1)/ \alpha}
((1-\alpha)/ \alpha)^{(1-\alpha)/ \alpha }}
\frac{\Gamma^{2}(\frac{1- \alpha }{2\alpha})}{\Gamma^{2}(\frac{1}{2\alpha})}
\right)^{\frac{\alpha }{(\alpha +1)}}.
\end{equation}
Here  $\jmath _{(n-3)/ 2}$ denotes the smallest positive zero of the Bessel function
$J_{(n-3)/2}$ and
$\Gamma $ is the Gamma function.
In particular,  since $\jmath_{-1/2}=\pi /2$, this limit is  $3 \pi ^{2} / 8$ for the parabola $\mathcal{P}$, for which  $A=1$, $n=2$ and $\alpha =1/2$.

The following is the first main result in this paper. For
simplicity of notation we will use $\tau _{\alpha}(Z)$ instead of
$\tau _{P_{\alpha}}(Z) $ to denote the first exit time of the
process $Z_{t}$ from $P_{\alpha}$.

\begin{theorem}\label{theorem1}
Let $ 0 < \alpha <1$ and let  $P_{\alpha}=\{ (x,Y)\in \RR{R} \times \RR{R}^{n-1}:
 x>0, |Y|<Ax^{\alpha} \}$. Then for
$z \in P_{\alpha}$,
\begin{eqnarray}
  \lim_{t\to\infty}t^{-(\frac{ 1- \alpha }{3+\alpha})}\log P_{z}[\tau _{\alpha}(Z)> t]
 = -(\frac{3+\alpha }{2+ 2 \alpha } )
(\frac{1+\alpha }{1-\alpha})^{(\frac{1- \alpha }{3+
\alpha} )} \pi ^{(\frac{2-2\alpha }{3+\alpha } )} l^{  (\frac{2+2\alpha }{3+\alpha }) } ,  \nonumber
\end{eqnarray}
where $l$ is the limit given by (\ref{brownianlimit}).
\end{theorem}

In the case of two dimensions and the parabola
$\mathcal{P}$, $A=1$
 we obtain by substituting $l=3\pi^{2} /8$ in Theorem \ref{theorem1}
\begin{eqnarray}
\lim_{t\rightarrow \infty} t^{-\frac{1}{7}}  \log \ P_{z}[\tau _{\mathcal{P}}(Z) >t] =
 - \frac{7 \pi ^{2}}{2^{25/ 7}}.\nonumber
\end{eqnarray}
Theorem \ref{theorem1} gives  the sharp order of integrability for iterated Brownian motion in these regions.

\begin{corollary}\label{corollary1}
Let $P_{\alpha}$ be a general parabolic domain as in Theorem \ref{theorem1}. Then for $z \in P_{\alpha}$,
$$ E_{z}\left[  \exp(b (\tau _{\alpha} (Z))^{\frac{  1- \alpha  }{\alpha +3}}) \right]$$
is finite if $$b<(\frac{3+\alpha }{2+ 2 \alpha } )
 (\frac{1+\alpha }{1-\alpha})^{ (\frac{1- \alpha
}{3+\alpha}) } \pi
^{(\frac{2-2\alpha }{3+\alpha } )} l^{  (\frac{2+2\alpha }{3+\alpha }) },$$
and it is infinite if
$$b>(\frac{3+\alpha }{2+ 2 \alpha } )
 (\frac{1+\alpha }{1-\alpha})^{ (\frac{1- \alpha
}{3+\alpha}) } \pi
^{(\frac{2-2\alpha }{3+\alpha } )} l^{  (\frac{2+2\alpha }{3+\alpha }) }.$$
\end{corollary}
 Note that we do not know what happens in the corollary when
$$
b=(\frac{3+\alpha }{2+ 2 \alpha } )
 (\frac{1+\alpha }{1-\alpha})^{ (\frac{1- \alpha
}{3+\alpha}) } \pi
^{(\frac{2-2\alpha }{3+\alpha } )} l^{  (\frac{2+2\alpha }{3+\alpha }) }.
$$
The proof of Corollary \ref{corollary1} follows from  Theorem \ref{theorem1} and the fact that
$$
E_{z}\left[  \exp(b (\tau _{\alpha} (Z))^{\frac{ 1- \alpha  }{\alpha +3}}) \right]=\int_{0}^{\infty}
(\frac{d}{dt}\exp (b t^{\frac{  1- \alpha }{\alpha +3}}))P_{z}[\tau _{\alpha} (Z)>t] dt.
$$
In \cite{DSmits}, DeBlassie and Smits studied the tail behavior of
the first exit time of the Brownian motion from twisted domains in
the plane. Let $D\subset \RR{R}^{2}$ be a domain whose boundary
consists of three curves (in polar coordinates)
\begin{eqnarray*}
\ & C_{1}:&   \theta =f_{1}(r), \ \  r\geq r_{1}\\
\ & C_{2}:&  \theta = f_{2}(r), \ \ r\geq r_{1} \\
\ & C_{3}:&   r=r_{1},  \ \ \ \ \ \  f_{2}(r)\leq \theta  \leq f_{1}(r)
\end{eqnarray*}
where $f_{1}$ and $f_{2}$ are smooth and the curves $C_{1}$ and $C_{2}$ do not cross:
$$
0<f_{1}(r)-f_{2}(r)<\pi, \ \ r\geq r_{1}.
$$
DeBlassie and Smits call $D$ a twisted domain if there
 are constants $r_{0}>0$, $\gamma >0$ and $p\in (0,1]$ and a smooth function $f(r)$ such
that the curves $f_{1}(r)$ and $f_{2}(r)$, $r\geq r_{0},$ are
obtained from $f(r)$  by moving out $\pm \gamma r^{p}$ units along
the normal to the curve $\theta=f(r)$ at the point whose polar
coordinates are $(r,f(r))$. They call $\gamma r^{p}$ the growth
radius and $\theta =f(r)$ the generating curve. DeBlassie and
Smits \cite[Theorem1.1]{DSmits} have the following tail behavior
of the first exit time of Brownian motion from twisted domains
$D\subset \RR{R}^{2}$ with growth radius $\gamma r^{p}$, $\gamma
>0, $ $0<p<1$
\begin{equation}\label{brownianlimit2}
\lim_{t\to \infty}t^{-(\frac{1-p}{1+p})}\log P_{z}[\tau_{D}>t]=-l_{1}=-\left[
\frac{\pi^{2p-1}}{\gamma 2^{2p}(1-p)^{2p}}
\right]^{\frac{2}{p+1}} C_{p}
\end{equation}
where
$$
C_{p}=(1+p)\left[
\frac{\pi^{2+p}}{8^{p}p^{2p}(1-p)^{1-p}}\frac{\Gamma ^{2p}
\left( \frac{1-p}{2p}\right)}{\Gamma ^{2p}\left( \frac{1}{2p}\right)}
\right]^{\frac{1}{p+1}}.$$
 For these domains  we  have the following theorem which is the last main result in this paper.

\begin{theorem}\label{theorem11}
Let $D\subset \RR{R}^{2}$ be a twisted domain with growth
 radius $\gamma r^{p}$, $\gamma >0,$ $0<p<1$.
Then $$
\lim_{t\to\infty} t^{-(\frac{1-p}{p+3})}\log P_{z}[\tau_{D}(Z)>t]=
-(\frac{3+p }{2+ 2 p } )
(\frac{1+p}{1-p})^{(\frac{1-p}{3+p})} \pi ^{(\frac{2-2p}{3+p } )}
l_{1}^{  (\frac{2+2p }{3+p }) } ,
$$
where $l_{1}$ is the limit given by (\ref{brownianlimit2}).
\end{theorem}
Notice the similarity of the limits in Theorems \ref{theorem1}
 and \ref{theorem11} except the constants $l$ and $l_{1}$.

To obtain our results we follow the general ideas of DeBlassie \cite{deblassie}
but with some key modifications.  We use
the asymptotics of $\frac{\partial }{\partial u} \frac{\partial }{\partial v}
P_{0}[\eta  _{(-u,v)}>t]$, where $\eta  _{(-u,v)}$ is the first exit time of  one-dimensional
Brownian motion
from the interval $(-u, v)$ and
integration by parts to get away from  assuming the asymptotics of the density of the exit times.

The paper is organized as follows. In $\S2$ we give some preliminaries needed in the proofs of
Theorems \ref{theorem1} and \ref{theorem11}.   In $\S3$, we prove
Theorems \ref{theorem1} and  \ref{theorem11}. In $\S4$ we derive several technical
 results that are used in the  asymptotics
of
$\frac{\partial }{\partial u} \frac{\partial }{\partial v}
P_{0}[\eta  _{(-u,v)}>t]$.

\section{Preliminaries}\label{prelim}
 In this section we prove some results which we will use in section \ref{main}.
In what follows  we will write $f \approx g$ and $f \lesssim g$ to mean
that
for some positive $C_{1}$ and $C_{2}$, $C_{1}\leq f/g \leq C_{2}$ and $f \leq C_{1} g$,
respectively. We will also write $f(t) \sim
g(t)$, as
$t\rightarrow \infty $,  to mean $f(t) / g(t) \rightarrow 1$, as $t\rightarrow \infty $.

\begin{lemma}\label{sumlemma}
Let $0< \beta \leq 1$.
Let $\xi$ be a positive random variable such that for some $c>0$, $-\log P[\xi >t] \sim ct^{\beta}$,
 as $t \rightarrow \infty$.
Then for independent copies $\xi _{1}$ and $\xi _{2}$ of $ \xi$, $-\log P[\xi_{1} +\xi _{2} >t] \sim
 ct^{\beta}$, as $t \rightarrow
\infty$.
\end{lemma}

\begin{proof} 
%By independence,
%\begin{eqnarray}
%P[\xi_{1} +\xi _{2} >t]&=& 1-\int_{0}^{t} P[\xi \leq t-y]d_{y}P[\xi \leq y]     \nonumber\\
%&=&1-\int_{0}^{t}(1- P[\xi > t-y])d_{y}P[\xi \leq y]    \nonumber\\
%&=&P[\xi >t] +\int_{0}^{t} P[\xi > t-y]d_{y}P[\xi \leq y]  \nonumber\\
%&\geq&P[\xi >t].    \nonumber
%\end{eqnarray}
The lower bound for $P[\xi_{1} +\xi _{2} >t]$ follows from the observation that $\xi_{1} +\xi _{2}$ is at least
$\xi_{1}$; so $P[\xi_{1} +\xi _{2} >t]\geq P[\xi >t].$

For an upper bound for $P[\xi_{1} +\xi _{2} >t]$, note that for any $\theta<c$,
$$E(\exp(\theta \xi ^{\beta}))< \infty .$$
Then by Chebyshev inequality and independence
\begin{eqnarray}
P[\xi_{1} +\xi _{2} >t]&\leq&e^{-\theta t^{\beta}}E(\exp (\theta ( \xi
_{1} +\xi _{2}) ^{\beta}))     \nonumber\\
&\leq&   e^{-\theta t^{\beta}}E(\exp (\theta ( \xi _{1}^{\beta} +\xi _{2} ^{\beta}))    \nonumber\\
&=& e^{-\theta t^{\beta}}E(\exp (\theta ( \xi  ^{\beta}))^{2} ,  \nonumber
\end{eqnarray}
where we have used the fact that
 for $0 < \beta \leq 1$ and $a,b $ positive real numbers,  $(a+b)^{\beta} \leq a^{\beta}+
b^{\beta}$. So,
\[\theta t^{\beta}-2\log E(\exp(\theta ( \xi  ^{\beta}))\leq  -\log P[\xi_{1} +\xi _{2} >t] \leq -\log P[\xi >t].
\]
Now divide by $ct^{\beta}$, let $t\rightarrow \infty $ and $\theta \uparrow c$
 to  get the desired conclusion.
\end{proof}

\begin{lemma}\label{sumlemma2}
Let $\beta>1$. Let $\xi$ be a positive random variable such that for some $c>0$, $-\log P[\xi >t] \sim ct^{\beta}$,
 as $t \rightarrow \infty$. Then
$$
 \lim_{t\rightarrow \infty} t^{-\beta}\log P[\xi_{1} +\xi _{2} >t]= -c 2^{1-\beta}.
$$
\end{lemma}

\begin{proof}
The proof of the lower bound for $P[\xi_{1} +\xi _{2} >t]$ follows from the fact that
$$
P[\xi_{1} +\xi _{2} >t] \geq (P[\xi > t/2])^{2}.
$$
Indeed, if $\xi_1$ and
   $\xi_2 $ are both at least $ t/2 $, then $\xi_1+\xi_2>t $.
   This implies that  for any $\epsilon >0$, $ P[\xi_{1} +\xi _{2} >t] \geq \exp (-2^{1-\beta}c t^\beta (1+\epsilon) )$ for $t$ large.

For the upper bound for $P[\xi_{1} +\xi _{2} >t]$ we use the Chebyshev inequality and the fact that for $\beta >1$ and $a,b$ positive real numbers,
$(a+b)^{\beta}\leq 2^{\beta -1}(a^{\beta}+b^{\beta})$.
\end{proof}

We use Lemma \ref{sumlemma} to derive the asymptotics of the
Laplace transform of $(\xi_{1} + \xi _{2})^{-2}$.
This is a special case of the  following theorem
(Kasahara \cite[Theorem 3]{kasahara} and Bingham,
Goldie and Teugels \cite[p. 254]{bgt}.)

\begin{theorem}[de Bruijn's Tauberian Theorem]

Let $X$ be a positive random variable such that for some positive $B_{1},\ B_{2}$ and $p$,
$$-B_{1}\leq \liminf_{x\rightarrow 0} x^{p}\log P[X \leq x] \leq \limsup_{x\rightarrow 0} x^{p}\log P[X \leq x]\leq -B_{2}  .$$
Then
\begin{eqnarray}
\  & \ &  -(p+1)(B_{1})^{1/(p+1)}p^{-p/(p+1)}
\leq \liminf_{\lambda \rightarrow \infty} \lambda^{-p/(p+1)}\log E e^{- \lambda X}\nonumber \\
\ & \ & \leq  \limsup_{\lambda \rightarrow \infty} \lambda^{-p/(p+1)}\log E e^{- \lambda X}
\leq -(p+1)(B_{2})^{1/(p+1)}p^{-p/(p+1)} .\nonumber
\end{eqnarray}

\end{theorem}

\begin{lemma}\label{lemma2}
Let $\xi$ be a positive random variable such that for some $c>0$, $-\log P[\xi >t] \sim ct^{\beta}$,
 as $t \rightarrow \infty$.
Let $\xi _{1}$ and $\xi _{2}$ be independent copies of $\xi$.  Then, for $0<\beta \leq 1$
\[
-\log E(\exp (-\frac{\lambda}{(\xi_{1}+\xi_{2})^{2}})) \sim  ((\beta +2) / 2)c^{2 / (2+\beta )}
(\beta / 2)^{- \beta / (\beta +2)}\lambda ^{\beta / (\beta +2)},
\]
as $  \lambda \rightarrow \infty,$ and  for $\beta>1$

$$-\log E(\exp (-\frac{\lambda}{(\xi_{1}+\xi_{2})^{2}})) \sim
  ((\beta +2) / 2)(c2^{1- \beta})^{2/(\beta +2)}(\beta / 2)^{- \beta / (\beta +2)}\lambda^{\beta/(\beta +2)}
 $$

as $\lambda \rightarrow \infty$.

\end{lemma}
\begin{proof} For $0<\beta \leq 1$,
\[
-\log P[\frac{1}{(\xi_{1}+\xi_{2})^{2}} \leq x]= -\log P[\xi_{1}+\xi_{2}>x^{-1 /2}]  \sim  cx^{-\beta / 2},
\ \  \mathrm{as} \ \ x\rightarrow 0,
\]
by Lemma \ref{sumlemma}. For $p= \beta / 2$ in de Bruijn's Tauberian Theorem we get
\[
-\log E(\exp(-\frac{\lambda}{(\xi_{1}+\xi_{2})^{2}})) \sim    ((\beta +2) / 2)c^{2 / (2+\beta )}
(\beta / 2)^{- \beta / (\beta +2)}\lambda ^{\beta / (\beta +2)},
\]
as $\lambda \rightarrow \infty.$

 For $\beta>1$, we use Lemma \ref{sumlemma2} and de Bruijn's Tauberian Theorem with,
$B_{1}=B_{2}=c/2^{\beta -1}$ and $p=\beta /2$.
\end{proof}

We also need the following application of de Bruijn's Tauberian Theorem.
\begin{lemma}\label{lemma3}
Let $X$ be a positive  random variable with density $$f(u)= \gamma u^{-2} e^{- \alpha / u^{\beta /2}},$$ then
$-\log \ P[X\leq x] \sim  \alpha x^{-\beta / 2}$,  as  $x\rightarrow 0$. In this case
$$
-\log \ E(e^{-\lambda X}) \sim ((\beta +2) / 2) \alpha ^{2 / (2+\beta )}
(\beta / 2)^{- \beta / (\beta +2)}\lambda ^{\beta / (\beta +2)},
$$
as $ \lambda \rightarrow \infty$.
\end{lemma}
\begin{proof}
\[
P[X\leq x]=\gamma \int_{0}^{x} u^{-2} e^{- \alpha / u^{\beta / 2}} du,
\]
and making the change of variables $v=u^{- \beta / 2 }$ gives $u=v^{- 2 / \beta }$ and  $du= -\frac{2}{\beta}
v^{- 2 / \beta-1} dv$ we get
\begin{eqnarray}
 P[X\leq x]&=&\gamma \frac{2}{\beta} \int_{x^{- \beta / 2}}^{\infty} v^{2 / \beta -1} e^{-\alpha v} dv \nonumber \\
&=&\gamma \frac{2}{\beta}  \frac{1}{\alpha} (\frac{1}{\alpha})^{1- 2 / \beta} \int_{\alpha
x^{- \beta / 2}}^{\infty}z^{2 / \beta  -1}e^{-z}dz \label{eqvar}\\
&=&\gamma \frac{2}{\beta} \alpha^{-\frac{2}{\beta}}(\alpha x^{- \beta / 2})^{2/ \beta  -1} e^{-\alpha  x^{- \beta / 2}}
[1+O( x^{ \beta / 2}) ]. \label{eqgrad}
\end{eqnarray}
In Equation (\ref{eqvar}) as $x\rightarrow 0$, $x^{- \beta / 2}\rightarrow \infty$.
Then equation (\ref{eqvar}) follows by changing variables $\alpha v=z$ and equation (\ref{eqgrad}) follows from the asymptotics
of incomplete gamma function as in Gradshteyn and Ryzhik \cite[p.942]{grad}.
Hence,
$$
x^{\frac{\beta}{2}} \log \ P[X\leq x]=x^{\frac{\beta}{2}}\ \log C_{1} +x^{\frac{\beta}{2}}C_{2}\ \log x -\alpha
+x^{\frac{\beta}{2}}\ \log (1+O(x^{\frac{\beta}{2}}))\ \ \rightarrow \ -\alpha ,
$$
as  $x\rightarrow 0$,
where $C_{1}=\gamma \frac{2}{\beta} \frac{1}{\alpha}$ and $C_{2}=- \frac{2-\beta}{2}$.
Now the desired conclusion follows from de Bruijn's Tauberian Theorem.
\end{proof}

\section{ Proof of main results }\label{main}

If $D\subset \RR{R} ^{n}$  is an open set, write
\[
\tau_{D}^{\pm}(z)=\inf \{ t\geq 0:\ \ X_{t}^{\pm} +z \notin D\},\]
and if $I\subset \RR{R}$ is an open interval, write
\[
\eta _{I}=\eta (I)= \inf \{ t\geq 0 :\ \  Y_{t}\notin I\}.
\]
Recall that $\tau _{D}(Z)$ stands for the first exit time of iterated Brownian motion from  $D$. As in DeBlassie \cite[\S3.]{deblassie}, we have by the continuity of the paths for $Z_{t}=z+X(Y_{t})$
\begin{eqnarray}
P_{z}[\tau _{D}(Z) > t]&=&P_{z}[Z_{s} \in D\  \mathrm{for}\ \ \mathrm{all} \  s\leq t]\nonumber \\
&=&P[z+X^{+}(0\vee Y_{s})\in D \  \mathrm{and}\  z+X^{-}(0\vee (-Y_{s}))\in D\ \nonumber \\
& &\ \mathrm{for}\ \ \mathrm{all}\  s\leq t ]\nonumber \\
&=&P[\tau _{D}^{+}(z) > 0\vee Y_{s} \   \mathrm{and} \   \tau _{D}^{-}(z) > 0\vee (-Y_{s})
 \ \mathrm{for}\ \ \mathrm{all}\ s\leq t]\nonumber \\
&=&P[-\tau _{D}^{-}(z) < Y_{s} < \tau _{D}^{+}(z) \ \mathrm{for}\ \ \mathrm{all}\ s\leq t]\nonumber \\
&=&P[\eta(-\tau _{D}^{-}(z) , \tau _{D}^{+}(z) ) >t],
\label{probeq}
\end{eqnarray}
and this equals for the parabola-shaped domains and twisted domains introduced above
\begin{equation}\label{eqddo}
=\int_{0}^{\infty} \!\int_{0}^{\infty}  \left( \frac{\partial }{\partial u} \frac{\partial }{\partial v}
 P_{0}[\eta _{(-u,v)}>t] \right) P[\tau_{D}(z) >u]P[\tau_{D}(z)  >v]dvdu.
\end{equation}
The equation (\ref{eqddo}) follows from Lemma \ref{lemmaA.3} below.

For the parabola-shaped domains $P_{\alpha} \subset \RR{R}^{n}$ we have
 $$
\lim_{t\to\infty}t^{-(\frac{1-\alpha } { 1+\alpha} )}
\log  P_{z}[\tau _{\alpha}> t]=- l,
$$
where $l$ is given in (\ref{brownianlimit}). Similarly for the twisted domains $D\subset \RR{R}^{2}$
$$
\lim_{t\to \infty}t^{-(\frac{1-p}{1+p})}\log P_{z}[\tau_{D}>t]=-l_{1},
$$
where $l_{1}$ is given in (\ref{brownianlimit2}).
Then our main theorem, Theorem \ref{theorem1}, follows from equation  (\ref{probeq}) and
 by substituting $(1-\alpha)/(1+\alpha)$ for $\beta$ and $l$ for $c$ in the following theorem. Similarly
Theorem \ref{theorem11} follows from the following theorem as well from equation (\ref{probeq})
and by substituting $(1-p)/(1+p)$ for $\beta$ and $l_{1}$ for $c$.

The following theorem, Theorem \ref{theorem31}, is more general than DeBlassie's theorem \cite[Theorem 4.4]{deblassie}, we use the asymptotics of the
distribution of the
random variables rather than the density of the random variables. For this, we use
integration by parts and the asymptotics of $ \frac{\partial }{\partial u} \frac{\partial }{\partial v}
P_{0}[\eta _{(-u,v)}>t]$.

\begin{theorem}\label{theorem31}
Let $0<\beta \leq 1$. Let $\xi$ be a positive random variable such that
$-\log P[\xi > t] \sim c t^{\beta}$, as $t\rightarrow \infty $. If $\xi _{1}$ and $\xi_{2}$ are independent copies of  $\xi$
and independent of the Brownian motion $Y,$ then
\[
-\log P[\eta_{(-\xi_{1},\xi_{2})} >t] \sim   (\frac{2+\beta}{2})c^{2/ (2+\beta)} \beta^{-  \beta / (2+ \beta) } \pi ^{2
\beta/ (2+\beta)} t^{\beta / (2+ \beta)},
\]
as $t \rightarrow  \infty$.
\end{theorem}

\begin{proof}

 We suppose $Y_{0}=x$ and
the probability associated with this will be $P_{x}$.
The distribution of $\eta _{(-u,v)}$ is well known:
\begin{equation}
P_{0}[\eta _{(-u,v)} >t]=\frac{4}{\pi }\sum_{n=0}^{\infty}\frac{1}{2n+1} \exp (-\frac{(2n+1)^{2}\pi ^{2}}{2(u+v)^{2}} t) \sin
\frac{(2n+1)\pi u}{u+v}, \label{distr}
\end{equation}
see Feller \cite[pp. 340-342]{feller}.\newline
Let $\epsilon >0$.
From Lemma \ref{lemmaA.1}, choose $M>0$ so large that
\begin{equation}
P_{x}[\eta _{(0,1)} >t] \approx   \frac{4}{\pi} e ^{- \frac{\pi ^{2} t}{2} } \sin \pi x, \   \mathrm{for} \  t\geq M \
\mathrm{uniformly}\  \ x \in (0,1). \label{papprox}
\end{equation}
Then choose $\delta <\frac{1}{2}$ so small that
\begin{equation}
\sin \pi x \approx x, \ \   x\in (0,\delta].\label{sinapprox}
\end{equation}
From the hypothesis choose $K>0$ so large that
\begin{eqnarray}
e^{-c(1+\epsilon)u^{\beta}}\leq P(\xi >u) \leq e^{-c(1-\epsilon) u^{\beta}} \ \ \mathrm{for}\ \ u\geq K. \label{dfapprox}
\end{eqnarray}
We further assume that $t$ is so large that $K<\delta \sqrt{t/M}$.

If $f$ is the probability density function of the random variable $\xi$,
by independence of $Y$, $\xi _{1}$ and $\xi _{2}$, and using scaling and translation invariance of Brownian motion
\begin{eqnarray}
P[\eta _{(-\xi _{1}, \xi _{2})}>t] & = & P_{0} [\eta _{(-\xi _{1}, \xi _{2})} > t ] \nonumber \\
&=&\int_{0}^{\infty} \!\int_{0}^{\infty} P_{0} [\eta _{(-u,v)}>t] f(u)f(v)dvdu \label{eqint}  \\
&=&\int_{0}^{\infty} \!\int_{0}^{\infty} P_{\frac{u}{(u+v)}}[\eta
_{(0,1)}>\frac{t}{(u+v)^{2}}]f(u)f(v)dvdu, \label{upperapp}
\end{eqnarray}
and this equals
\begin{eqnarray}
=\int_{0}^{\infty} \!\int_{0}^{\infty} \left( \frac{\partial }{\partial u} \frac{\partial }{\partial v}
 P_{0}[\eta _{(-u,v)}>t]\right) P[\xi >u]P[\xi >v]dvdu.
\label{eqdd}
\end{eqnarray}
Equation (\ref{eqdd}) follows from Lemma \ref{lemmaA.3} using integration by parts in (\ref{eqint}).
By Lemma \ref{lowerintegral} the integral in equation (\ref{eqdd}) over the set $ u+v\geq \sqrt{t/M}$ satisfies
\begin{eqnarray}
\int \int_{u+v\geq \sqrt{t/M}} \left( \frac{\partial }{\partial u} \frac{\partial }{\partial v} P_{0} [\eta _{(-u,v)}>t]\right) P[\xi >u]P[\xi >v]dvdu
  \nonumber \\
\gtrsim  -e^{-C_{0}(\sqrt{t/M})^{\beta}}.\label{goodlower1}
\end{eqnarray}
For large $t$, let
\[
A=\left\{ (u,v):\ \ K\leq u\leq \delta \sqrt{\frac{t}{M}} ,\ \frac{1-\delta}{\delta} u \leq v\leq \sqrt{\frac{t}{M}}-u \right\}.
\]
For $ u+v\leq \sqrt{t/M}$ we have the lower bound approximation from Lemma \ref{lemmaA.2} which says
for $ u+v\leq \sqrt{t/M}$
\begin{equation}
 \frac{\partial }{\partial u} \frac{\partial }{\partial v} P_{0}[\eta _{(-u,v)} > t]\gtrsim\
\exp(-\frac{\pi ^{2} t}{2 (u+v)^{2}}) \frac{uv}{(u+v)^{4}}\sin \frac{\pi
u}{u+v}. \label{lowerbound}
\end{equation}

 Now $A\subset \{(u,v): u+v\leq \sqrt{t/M} \}$.
Then by equations (\ref{eqdd})-(\ref{lowerbound})  we get
\begin{eqnarray}
P[\eta _{(-\xi _{1}, \xi _{2})}>t]&\gtrsim &\int \int_{A} \left( \frac{\partial }{\partial u} \frac{\partial }{\partial v} P_{0}[\eta _{(-u,v)}>t]\right)  P[\xi >u]P[\xi >v]dvdu\nonumber\\
& + & -e^{-C_{0}(\sqrt{t/M})^{\beta}}.\nonumber
\end{eqnarray}
We will show below that the integral on the right hand side of the last inequality is $\gtrsim \exp (-C_{1} t^{\beta /(\beta +2)})$ for some
$C_{1}$ positive. Thus, using the fact that for $t $ large $$ \exp (-C_{1} t^{\beta /(\beta +2)})-\exp(-C_{0}(\sqrt{t/M})^{\beta})
\gtrsim \exp (-C_{1} t^{\beta /(\beta +2)}),$$ in the
rest of the proof we will omit the second term, $\exp(-C_{0}(\sqrt{t/M})^{\beta})$, in the above inequality in finding the lower bound estimate for
$P[\eta _{(-\xi _{1}, \xi _{2})}>t]  .$

On the set $A$, since $\delta < 1/2 $, we have $v\geq (\frac{1}{\delta} -1)u>u>K$ and $u+v >\frac{u}{\delta}$, this gives
$\frac{u}{u+v}\leq \delta. $ By equations  (\ref{sinapprox}), (\ref{dfapprox}) and (\ref{lowerbound}), $P[\eta _{(-\xi _{1}, \xi _{2})}>t]$ is
\begin{eqnarray}
\gtrsim  \int_{K}^{\delta \sqrt{t/ M}} \int_{(1-\delta)u/
\delta}^{\sqrt{t/ M}-u} \exp (-\frac{\pi ^{2} t}{2 (u+v)^{2}}) \frac{u^{2}v }{ (u+v)^{5}}
e^{-(1+\epsilon)c(u^{\beta} +v^{\beta})} dvdu, \nonumber
\end{eqnarray}
and this is
\begin{equation}
\gtrsim  \int_{K}^{\delta \sqrt{t/ M}}
\int_{(1-\delta) u/ \delta}^{\sqrt{t/ M}-u} \exp(-\frac{\pi ^{2} t}{2 (u+v)^{2}}) \frac{u}{(u+v)^{5}}
e^{-(1+\epsilon)c(1+\delta ^{\beta})(u+v)^{\beta}}dvdu. \label{deltk}
\end{equation}
Inequality  (\ref{deltk}) follows from the fact that $u< \delta (u+v)$ over the set $A$ which gives $u^{\beta} +v^{\beta} \leq (1+\delta
^{\beta})(u+v)^{\beta},$ and $uv\geq K^{2}$ over $A$.
Changing the variables $x=u+v, z=u$ the integral is
\[
\approx \int_{K}^{\delta \sqrt{t/ M}} \int_{z/ \delta}^{\sqrt{t/ M}}\exp(-\frac{\pi ^{2} t}{2 x^{2}}) \frac{z}{x^{5}}
e^{-(1+\epsilon)c(1+\delta ^{\beta}) x^{\beta}}dxdz,
\]
and reversing the order of integration
\begin{eqnarray}
\ &\ & =\int_{K/ \delta}^{\sqrt{t/ M}} \int_{K}^{\delta x}\ \frac{z}{ x^{5} } \exp(-\frac{\pi ^{2} t}{2 x^{2}})
e^{-(1+\epsilon)c(1+\delta ^{\beta}) x ^{\beta}}dzdx \nonumber \\
 \ &\ &\approx  \int_{K/ \delta}^{\sqrt{t/ M}} \frac{1}{x^{5}}  \exp (-\frac{\pi ^{2} t}{2 x^{2}})
e^{-(1+\epsilon)c(1+\delta ^{\beta}) x^{\beta}}\  (\delta ^{2} x^{2}-K^{2})dx \nonumber \\
\ &\ & \geq \int_{ 2 K/ \delta}^{\sqrt{t/ M}}\frac{1}{x^{5}} \exp (-\frac{\pi ^{2} t}{2 x^{2}})
e^{-(1+\epsilon)c(1+\delta ^{\beta}) x^{\beta}}\  (\delta ^{2} x^{2}-K^{2})dx , \  \  t \   \ \mathrm{large} \nonumber  \\
\ & \ & \geq  \int_{ 2 K/ \delta}^{\sqrt{t/ M}}\frac{1}{x^{5}} \exp (-\frac{\pi ^{2} t}{2 x^{2}})
e^{-(1+\epsilon)c(1+\delta ^{\beta}) x ^{\beta}}\  (\delta ^{2} x^{2}-\frac{\delta ^{2} x^{2}}{4})dx. \nonumber
\end{eqnarray}
Changing variables $u=x^{-2}$ this is
\begin{eqnarray}
\ &\ &\approx \int_{ M/ t}^{\delta ^{2} /4K^{2}}  \exp(-\pi ^{2} tu  / 2 )
e^{-(1+\epsilon)c(1+\delta ^{\beta})(u)^{- \beta /2}}\ du  \nonumber  \\
\ &\ &\geq (M/t)^{2}\int_{ M/ t}^{\delta ^{2} /4K^{2}} u^{-2} \exp(-\pi ^{2} tu  / 2 )
e^{-(1+\epsilon)c(1+\delta ^{\beta})(u)^{- \beta /2}} du .    \nonumber
\end{eqnarray}
Thus we have for $t$ large,
\begin{eqnarray}
\ &\ & P[\eta _{(-\xi _{1}, \xi _{2})}>t] \nonumber\\
\ &\ & \gtrsim
(M/t)^{2}\int_{ M/ t}^{\delta ^{2} /4K^{2}} u^{-2} \exp(-\pi ^{2} tu  / 2 )
e^{-(1+\epsilon)c(1+\delta ^{\beta})(u)^{- \beta /2}}\ du. \label{eqlast}
\end{eqnarray}
We can disregard $(M/t)^{2}$ since we will take $\log $ and divide by $t^{\beta /( 2+ \beta) }$
and let $t\rightarrow \infty$.  That is,
$t^{-\beta/(\beta +2)} \log t \rightarrow 0$, as $t\rightarrow \infty$.
Now for some $c_{1}>0$
\begin{eqnarray}
\ & \ & \int_{\delta ^{2} /4K^{2}}^{\infty} u^{-2} \exp(-\pi ^{2} tu  / 2 )
e^{-(1+\epsilon)c(1+\delta ^{\beta})u^{- \beta /2}}\ du   \nonumber\\
\ &\ & \leq  e^{- \pi^2 \delta ^{2}t /8K^{2}}\int_{\delta ^{2}
/4K^{2}}^{\infty} u^{-2} e^{-(1+\epsilon)c(1+\delta ^{\beta})u^{- \beta /2}}\ du \nonumber\\
\ &\ & \lesssim  e^{-c_{1} t}. \label{last1}
\end{eqnarray}
Changing variables $u=v^{-2}$
\begin{eqnarray}
\ & \ &\int_{0}^{M/t} u^{-2} \exp(-\pi ^{2} tu  / 2 )
e^{-(1+\epsilon)c(1+\delta ^{\beta}) u^{- \beta /2}}\ du  \nonumber \\
 \ &\ &\leq   \int_{0}^{M/t} u^{-2} e^{-(1+\epsilon)c(1+\delta ^{\beta}) u^{- \beta /2}}\ du  \nonumber \\
  \ &\ & =  2\int_{\sqrt{t/M}}^{\infty} v e^{-(1+\epsilon)c(1+\delta ^{\beta}) v^{\beta}} dv \nonumber  \\
 \ &\ & \lesssim  (t/M)^{\beta /2  (2 / \beta -1)} e^{-(1+\epsilon)c(1+\delta ^{\beta})(t/M)^{\beta /2}}. \label{last2}
\end{eqnarray}
The inequality (\ref{last2}) follows from the asymptotics of the incomplete gamma function as in  Gradshteyn and Ryzhik  \cite[p.942]{grad}.
By  Lemma \ref{lemma3}, with  $\alpha =c(1+\epsilon)(1+\delta ^{\beta})$ and $\lambda =\frac{\pi ^{2}t}{2}$
\begin{eqnarray}
 \ & \ &  -\log \ \int_{0}^{\infty} u^{-2} \exp(-\pi ^{2} tu  / 2 )
e^{-(1+\epsilon)c(1+\delta ^{\beta}) u^{- \beta /2}}\ du   \ \nonumber \\
  \ &\ & \sim \csti{\frac{\beta}{2}}{c(1+\epsilon)(1+\delta ^{\beta})}{\frac{\pi ^{2}t}{2}}   \nonumber \\
\ &\ &= (\frac{2+\beta }{ 2})(c(1+\epsilon)(1+\delta ^{\beta}))^{\frac{2}{
(2+\beta )}}\beta ^{- \frac{\beta }{ (2+ \beta )} }
\pi ^{\frac{2 \beta}{ (2+\beta )}} t^{\frac{\beta } { (2+ \beta)}}, \label{logapp}
\end{eqnarray}
as $t \rightarrow \infty$.

By  equations (\ref{goodlower1}), (\ref{last1})- (\ref{logapp})
and from the fact that $\beta  / (\beta +2) < \beta /2$ and $\beta
/ (\beta +2) <1$, we can rewrite equation (\ref{eqlast}) for large
$t$ as $ P[\eta _{(-\xi _{1}, \xi _{2})}>t]$ is
\begin{equation}\label{upapp}
 \gtrsim \exp \left( -(1+ \epsilon)(\frac{2+\beta }{
2})(c(1+\epsilon)(1+\delta ^{\beta}))^{\frac{2}{ (2+\beta )}}\beta
^{- \frac{\beta }{ (2+ \beta )} } \pi ^{\frac{2 \beta}{ (2+\beta
)}} t^{\frac{\beta } { (2+ \beta)}} \right).
\end{equation}
Now we give an upper bound.
By equations (\ref{papprox})  and (\ref{upperapp})
\begin{eqnarray}
P[\eta _{(-\xi _{1}, \xi _{2})}>t] & = &
\ \ \int_{0}^{\infty}  \! \int_{0}^{\infty}
 P_{u / (u+v)} [\eta _{(0,1)}>\frac{t}{(u+v)^{2}}]\ f(u)f(v)dvdu  \nonumber \\
&\lesssim & \int \! \! \int _{u+v \leq \sqrt{t/ M}} e^{-\frac{\pi ^{2} t}{2 (u+v)^{2}}}\ f(u)f(v)dvdu \nonumber \\
&+ & \int \! \! \int _{u+v \geq  \sqrt{t/ M}}\ \ \ f(u)f(v)dvdu \nonumber \\
& \leq  & E \left[ \exp(-\frac{\pi ^{2} t}{2(\xi_{1}+\xi_{2})^{2}}) \right]+P(\xi_{1}+\xi_{2} \geq  \sqrt{t/ M}). \label{lowerapp}
\end{eqnarray}
By Lemma \ref{sumlemma} our random variables $\xi_{1}$ and $\xi_{2}$ satisfy
\begin{eqnarray}
-\log P[\xi_{1} + \xi_{2} >t]  \sim   ct^{\beta},\ \  \mathrm{as}\ \  t \rightarrow \infty .\nonumber
\end{eqnarray}
Then by Lemma \ref{lemma2} for $t$ large, equation (\ref{lowerapp}) becomes
\begin{eqnarray}
 \ &\ &P[\eta _{(-\xi _{1}, \xi _{2})}>t]\nonumber\\
\ &\ & \lesssim \exp\left( -(1- \epsilon)(\frac{2+\beta }{ 2})c^{\frac{2}{
(2+\beta )}}\beta ^{- \frac{\beta }{ (2+ \beta )} }
\pi ^{\frac{2 \beta}{ (2+\beta )}} t^{\frac{\beta } { (2+ \beta)}}\right) \nonumber \\
 \ &\ & + \ \exp \left( -c(1-\epsilon)(t/ M)^{\beta /2}\right) \nonumber \\
 \ &\ &\lesssim \exp \left( -(1- \epsilon)(\frac{2+\beta }{ 2}) c^{\frac{2}{
(2+\beta )}}\beta ^{- \frac{\beta }{ (2+ \beta )} }
\pi ^{\frac{2 \beta}{ (2+\beta )}} t^{\frac{\beta } { (2+ \beta)}} \right) .\nonumber
\end{eqnarray}
Combined with inequality (\ref{upapp}), this gives
\begin{eqnarray}
 \ &\ &-(1+ \epsilon)(\frac{2+\beta }{ 2})(c(1+\epsilon)(1+\delta ^{\beta}))^{\frac{2}{
(2+\beta )}} \beta^{- \frac{\beta }{ (2+ \beta )} }
\pi ^{\frac{2 \beta}{ (2+\beta )}}  \nonumber \\
 \ &\ & \leq   \liminf_{t \rightarrow \infty}   \   t^{ -\beta / (2+ \beta)}\  \log\ P[\eta _{(-\xi _{1}, \xi _{2})}>t]\nonumber\\
 \ &\ &\leq   \limsup_{t \rightarrow \infty}  \ t^{ -\beta / (2+ \beta)}\  \log\ P[\eta _{(-\xi _{1}, \xi _{2})}>t] \nonumber \\
 \ &\ &\leq -(1-\epsilon)(\frac{2+\beta }{ 2}) c^{\frac{2}{ (2+\beta )}} \beta^{- \frac{\beta }{ (2+ \beta )} }
\pi ^{\frac{2 \beta}{ (2+\beta )}}. \nonumber
\end{eqnarray}
Let $\epsilon \rightarrow 0$ and $\delta \rightarrow 0$  to get the desired conclusion.
\end{proof}

As a corollary of  Theorem \ref{theorem31} we obtain  DeBlassie's result  \cite[Theorem 4.4]{deblassie}.

\begin{corollary}
Let $0<\beta \leq 1$. Let $\xi$ be a positive random variable with a density function $f$ such that
$-\log f(t)  \sim ct^{\beta }$, as $t\rightarrow \infty$. If $\xi _{1}$ and $\xi_{2}$ are independent copies of  $\xi$
and independent of the Brownian motion $Y,$ then
\[
-\log P[\eta_{(-\xi_{1},\xi_{2})} >t] \sim   (\frac{2+\beta}{2})c^{2/ (2+\beta)} \beta^{-  \beta / (2+ \beta) } \pi ^{2
\beta/ (2+\beta)} t^{\beta / (2+ \beta)},
\]
as $t \rightarrow  \infty$.
\end{corollary}
The proof is immediate after we observe that if $-\log f(t)  \sim ct^{\beta }, $ as $t \rightarrow \infty$,
 then $$ -\log P[\xi > t] \sim c t^{\beta},\  \mathrm{as}\  t\rightarrow \infty. $$
We have the following corollary of Theorem \ref{theorem31}.
\begin{theorem}\label{theorem32}
Let $0<\beta \leq 1$. Let $D$ be a domain in $\RR{R}^{n}$. Let $\tau_{D}$ denote the first exit time of the Brownian motion from $D$ and satisfy
for $z\in D$,
$\lim_{t\to\infty}t^{-\beta}\log P_{z}[\tau _{D}> t]=-  c $ for some $c$ positive. If $\tau_{D}(Z)$ denotes the first exit time of the iterated Brownian motion from $D$, then
\[
\lim_{t\to\infty}t^{-\beta / (2+ \beta)}\log P[\tau_{D}(Z) >t]  = -(\frac{2+\beta}{2})c^{2/ (2+\beta)} \beta^{-  \beta / (2+ \beta) } \pi ^{2
\beta/ (2+\beta)}.
\]

\end{theorem}
Actually, Theorem \ref{theorem11} follows from  Theorem \ref{theorem32} above by substituting $(1-p)/(1+p)$ for $\beta$ and for $c$ the
limit in equation (\ref{brownianlimit2}). Theorem \ref{theorem1} as well follows from Theorem \ref{theorem32} in a similar way.

\begin{theorem}
Let $\beta >1$. Let $\xi$ be a positive random variable  such
that $-\log P[\xi > t] \sim c t^{\beta}$, as $t\rightarrow \infty $. If $\xi _{1}$ and $\xi_{2}$ are independent copies of  $\xi$
and independent of the Brownian motion $Y,$ then
\begin{eqnarray}
 \ &\ &-(\frac{2+\beta }{ 2}) c^{\frac{2}{
(2+\beta )}} \beta^{- \frac{\beta }{ (2+ \beta )} }
\pi ^{\frac{2 \beta}{ (2+\beta )}} \nonumber \\
\ & \ &  \leq   \liminf_{t \rightarrow \infty}   \   t^{ -\beta / (2+ \beta)}\  \log\ P[\eta _{(-\xi _{1}, \xi _{2})}>t]\nonumber\\
 \ &\ &\leq   \limsup_{t \rightarrow \infty}  \ t^{ -\beta / (2+ \beta)}\  \log\ P[\eta _{(-\xi _{1}, \xi _{2})}>t] \nonumber \\
 \ &\ &\leq -(\frac{2+\beta }{ 2})(c2^{1-\beta})^{\frac{2}{ (2+\beta )}} \beta^{- \frac{\beta }{ (2+ \beta )} }
\pi ^{\frac{2 \beta}{ (2+\beta )}}. \nonumber
\end{eqnarray}
\end{theorem}

\begin{proof} The result follows
from Lemma \ref{sumlemma2} , de Bruijn's Tauberian
Theorem and the proof of the Theorem \ref{theorem31},
using the well-known fact that for $\beta >1$ and
$a,\ b$ positive real numbers, $(a+b)^{\beta}\leq 2^{\beta -1}(a^{\beta}+b^{\beta}) $.
\end{proof}
 Comparing with a ball inside any domain $D$, we see that the analogue of the Theorem \ref{theorem32} does not hold for $\beta >1$. Let $B\subset D$ be a ball centered at $x$, then for the Brownian motion $X_{t}$ started at $x\in D$, 
 $-\log P_{x}[\tau_{D} > t]\lesssim \lambda_{B} t$ for large $t$ where $\lambda_{B}$ is the first eigenvalue of the Dirichlet Laplacian for $B$. This implies that the statement ` $ -\log P_{x}[\tau_{D} > t]\sim ct^{\beta}$, as
$t\rightarrow \infty$ '
 cannot be  true for $\beta >1$.

\section{Asymptotics}\label{approximation}

In this Section we will prove some lemmas that were used in section \ref{main}.
The following lemma is proved in \cite[Lemma A1]{deblassie} (it also follows from more general results on
``intrinsic ultracontractivity").  We include it for completeness.%\newline

\begin{lemma}\label{lemmaA.1}
As $t\rightarrow \infty$,
\[
P_{x}[\eta _{(0,1)} >t]\  \sim \  \frac{4}{\pi}e^{-\frac{\pi ^{2} t}{2} } \sin \pi x,
\ \  \mathrm{uniformly} \ \ \mathrm{for}\ \  x\in (0,1).
\]
\end{lemma}
We will need asymptotics of $ \frac{\partial }{\partial u} \frac{\partial }{\partial v} P_{0}(\eta _{(-u,v)} > t)$  for $(u,v)\in A$,
where we define $A$ for $\delta < 1 /2 ,\  K >0 $ and $ M>0$ as
\[
A=\left\{(u,v):\ \ K\leq u\leq \delta \sqrt{\frac{t}{M}} ,\ \frac{1-\delta}{\delta} u \leq v\leq \sqrt{\frac{t}{M}}-u \right\}.
\]
\begin{lemma}\label{lemmaA.2}
Let $B=\{ (u,v): \ t/ (u + v)^{2} >M\}$ for $M$ large. On $B $ we have
\begin{eqnarray}
 \ &\ & \frac{\partial }{\partial u} \frac{\partial }{\partial v} P[\eta _{(-u,v)} > t]\nonumber\\
  \ &\ & \approx \ \ 4 \exp(-\frac{\pi ^{2} t}{2 (u+v)^{2}})\left( (\sin \frac{\pi
u}{u+v}) (\frac{1
}{(u+v)^{4}})(\frac{\pi^{3}t^{2}}{(u+v)^{2}}
-3\pi t  \right. \nonumber \\
 \ &\ &\left. +\pi uv) +(\cos \frac{\pi u}{u+v})(\frac{1}{(u+v)^{3}}) (\frac{\pi^{2}t (v-u)}{(u+v)^{2}}+u-v) \right). \nonumber
\end{eqnarray}

Moreover, on $B$
\begin{equation}
 \frac{\partial }{\partial u} \frac{\partial }{\partial v} P_{0}[\eta _{(-u,v)} > t]\gtrsim\
\exp(-\frac{\pi ^{2} t}{2 (u+v)^{2}}) \frac{uv}{(u+v)^{4}}\sin \frac{\pi
u}{u+v}. \label{lowerdiff}
\end{equation}
\end{lemma}
\begin{proof}
If we  take the derivative of
\begin{eqnarray}
P_{0}[\eta _{(-u,v)} >t]=\frac{4}{\pi }\sum_{n=0}^{\infty}\frac{1}{2n+1} \exp (-\frac{(2n+1)^{2}\pi ^{2}}{2(u+v)^{2}} t) \ \sin \
\frac{(2n+1)\pi u}{u+v}, \nonumber
\end{eqnarray}
term by term w.r.t $u$, we get
\begin{eqnarray}
\ &\ & \frac{\partial }{\partial u}  P_{0}[\eta _{(-u,v)} >t]\label{duu}\\
\ &\ & = \frac{4}{\pi }\sum_{n=0}^{\infty} \  \exp (-\frac{(2n+1)^{2}\pi ^{2}}{2(u+v)^{2}} t)
\left[ \frac{(2n+1)\pi ^{2}t}{(u+v)^{3}} \ \sin \
\frac{(2n+1)\pi u}{u+v}\right. \nonumber\\
 \ &\ & \left. \ \ \ \ \ \ \ \ \ \ \ +   \frac{\pi v}{(u+v)^{2}}\ \cos\ \frac{(2n+1)\pi u}{u+v} \right] ,\nonumber
\end{eqnarray}
and if we take derivative of equation (\ref{duu}) with respect to $v$, then
\begin{eqnarray}
\ &\ & \frac{\partial }{\partial u} \frac{\partial }{\partial v} P_{0}[\eta _{(-u,v)} >t]\nonumber\\
\ &\ & = 4 \ \sum_{n=0}^{\infty} \  \exp (-\frac{(2n+1)^{2}\pi ^{2}}{2(u+v)^{2}} t)
\left\{\left[  \ \sin \ \frac{(2n+1)\pi u}{u+v}\ (\frac{1}{(u+v)^{4}}) \right. \right. \nonumber \\
 \ &\ & \left. \left.\ \ \ \ \ \ \  \times \ (\frac{\pi ^{3}(2n+1)^{3} t^{2}}{(u+v)^{2}} - \ 3\pi (2n+1)t \ + \ (2n+1)\pi uv)\ \
\right] \right.\nonumber\\
 \ &\ & +  \left. \left[ \ \cos\ \frac{(2n+1)\pi u}{u+v}\ (\frac{1}{(u+v)^{3}})(\frac{(2n+1)^{2}\pi ^{2}t (v-u)}{(u+v)^{2}} +\ \ u\ -\ v) \right]
\right\}. \nonumber
\end{eqnarray}
The term by term differentiation is admissible since we have the exponential function in each term.

The asymptotics in the lemma related to $\sin$ terms follow from the proof of Lemma \ref{lemmaA.1} in Deblassie 
\cite[Lemma A1]{deblassie}, which says that for some $c>0$
$$ \left|  \sin (2n+1)\pi x /  \sin \pi x  \right|\leq c(2n+1)^2,  \ \mathrm{ uniformly\   for\  all}\  x\in (0,1).$$
 %and that the fractions are bounded on $A$
The asymptotics for the
cosine terms
follow from the fact that by induction on $n$,
 $$ \left|  \cos (2n+1)\pi x /  \cos \pi x  \right|\leq 2n+1,  \ \mathrm{ uniformly\   for\  all}\  x\in (0,1).$$

The lower bound asymptotics follow from  the fact that for  $u+v <\sqrt{t /M}$
$$\frac{\pi^{3}t^{2}}{(u+v)^{2}}
-3\pi t \geq 0 , $$
 $$(\frac{\pi^{2}t }{(u+v)^{2}}-1)) >0,$$  and

$$ (v-u)\cos \frac{\pi u}{u+v} \geq 0 \ \ \mathrm{for} \ \ \mathrm{all}  \ \ u,v >0.$$
The last statement is valid since for  $u \leq v $ we have $(v-u)\geq 0$ and $u/(u+v)\leq 1/2$ which gives  $\cos \frac{\pi u}{u+v}\geq 0$, hence their product is positive.
For $u\geq v$ we have $v-u \leq 0$ and $\cos \frac{\pi u}{u+v} \leq 0$, so their product is positive. Hence we get the lower bound
asymptotics for $t/(u+v)^{2}$  large.
\end{proof}
Note that we obtain lower bound asymptotics of $\frac{\partial }{\partial u} \frac{\partial }{\partial v} P_{0} [\eta _{(-u,v)}>t]$ on the set $A$, because $A\subset B.$
\begin{lemma}\label{lowerintegral}
Let $\beta >0$. Let $\xi$ be a positive random variable such
that $-\log P[\xi > t] \sim c t^{\beta}$, as $t\rightarrow \infty $. If $\xi _{1}$ and $\xi_{2}$ are independent copies of  $\xi$
and independent of the Brownian motion $Y$, then
for some $ C_{0}>0 $
\begin{eqnarray}
\ &\ & \int \int_{u+v\geq \sqrt{t/M}} \left( \frac{\partial }{\partial u} \frac{\partial }{\partial v} P_{0} [\eta _{(-u,v)}>t]\right)  P[\xi >u]P[\xi >v]dvdu \nonumber\\
\ & \  &\gtrsim -e^{-C_{0}(\sqrt{t/M})^{\beta}} .\nonumber
\end{eqnarray}
\end{lemma}

\begin{proof}
%If we look at the formula for $\frac{d}{du} \frac{d}{dv} P_{0} (\eta _{(-u,v)}>t)$ in Lemma \ref{lemmaA.2}
We divide the set $u+v \geq \sqrt{t/M} $ into the following subsets
$$ A_{k}=\{ (u,v) : k\sqrt{t/M} \leq u+v < (k+1)\sqrt{t/M})  \}.$$
Now
\begin{eqnarray}
\! &\! & \int \! \int_{u+v\geq \sqrt{t/M}} \left( \frac{\partial }{\partial u} \frac{\partial }{\partial v} P_{0} [\eta _{(-u,v)}>t] \right) \ P[\xi >u]P[\xi >v]dvdu \nonumber \\
\!  &\! &=\sum_{k=1}^{\infty}
\int\!  \int _{A_{k}} \left( \frac{\partial }{\partial u} \frac{\partial }{\partial v} P_{0} [\eta _{(-u,v)} >t ]\right) \ P[\xi >u]P[\xi >v]dvdu \nonumber \\
\!  &\! &\gtrsim - \sum_{k=1}^{\infty}\sum_{n=0}^{\infty}
\int \! \int _{A_{k}}  (U +V) P[\xi >u]P[\xi >v]dvdu, \label{UVintegral}
\end{eqnarray}
where
\begin{eqnarray}
 U& = & \exp (-\frac{(2n+1)^{2}\pi ^{2}}{2(u+v)^{2}} t) \ (\frac{1}{(u+v)^{4}})\nonumber \\
 \ & \ & \times  \left| (\frac{\pi ^{3}(2n+1)^{3} t^{2}}{(u+v)^{2}}\
- \ 3\pi (2n+1)t + \ (2n+1)\pi uv ) \right|, \nonumber
\end{eqnarray}
and
\begin{eqnarray}
\ &\  & V =  \exp (-\frac{(2n+1)^{2}\pi ^{2}}{2(u+v)^{2}} t) \left|  \ (\frac{1}{(u+v)^{3}})(u- v
 +   \frac{(2n+1)^{2}\pi ^{2}t (v-u)}{(u+v)^{2}} ) \right|.\nonumber
\end{eqnarray}
We work on each term in equation (\ref{UVintegral}).
On the set $A_{k}$,
$$  \exp (-\frac{(2n+1)^{2}\pi ^{2}}{2(u+v)^{2}} t)\leq e^{-((2n+1)/(k+1))^{2} M/2}, $$
and terms like the reciprocal powers of $u+v$ are less than  the reciprocal of $k$ with the same
power. Hence we  have from equation (\ref{UVintegral}) for $t$ large
\begin{eqnarray}
\ &\ & \sum_{k=1}^{\infty}\sum_{n=0}^{\infty}
\int \int _{A_{k}}  (U +V) P[\xi >u]P[\xi >v]dvdu \nonumber \\
\ &\ & \lesssim \sum_{k=1}^{\infty}\sum_{n=0}^{\infty}
\int \int _{A_{k}} e^{-((2n+1)/(k+1))^{2} M/2}  \nonumber\\
  &\ &\times  \left[ \left( (\frac{1}{k^{4}} )(\frac{\pi ^{3}(2n+1)^{3} t^{2}}{k^{2}} +
 3\pi (2n+1)t)+ ((2n+1)\pi \frac{1}{k^{2}})\right) \right. \nonumber \\
\ &\ & \ \ \ \ \ \ \  \left.    + \left(
\frac{1}{(k)^{2}}(\frac{(2n+1)^{2}\pi ^{2}t }{k^{2}}\  \ + \  1)
\right) \right] P[\xi >u]P[\xi >v] dvdu .\label{lasteqn}
\end{eqnarray}
Now $A_{k} \subset A_{k}^{1} \cup A_{k}^{2}$, where
$$A_{k}^{1}=\{ (u,v):k/2 \sqrt{t/ M} \leq u \leq (k+1) \sqrt{t/ M} \  \mathrm{and} \ 0\leq v \leq (k+1)/2 \sqrt{t/ M}\} $$
$$A_{k}^{2}= \{
(u,v):k/2 \sqrt{t/ M} \leq v \leq (k+1) \sqrt{t/ M} \   \mathrm{and} \  0\leq u \leq (k+1)/2 \sqrt{t/ M}\}.
$$
From the asymptotics of $P[\xi >t]$, for $t$ large there is some $C_{0}>0$ independent of $k$ such that 
%on the sets containing $A_{k}$ on which either $u$, or $v$ is large, we have for some 
$P[\xi >u] P[\xi >v] \leq \exp (-C_{0}u^{\beta}) $, for $(u,v) \in A_{k}^{1}$. Thus by symmetry, same holds for  $(u,v) \in A_{k}^{2}$ 
\begin{eqnarray}
\ &\ &\int \int _{A_{k}}P[\xi >u]P[\xi >v]
\leq  \int \int _{A_{k}^{1} \cup A_{k}^{2}}P[\xi >u]P[\xi >v]\nonumber \\
 \ &\ & \lesssim
 (k+1)/2 \sqrt{t/ M}\int_{k/2 \sqrt{t/ M}}^{(k+1) \sqrt{t/ M}} e^{-C_{0}u^{\beta}}du \nonumber \\
\ &\ & \lesssim (k+1)/2 \sqrt{t/ M}(k+2)/2 \sqrt{t/ M}e^{-C_{0} (\frac{k}{2} \sqrt{t/M})^{\beta}},\nonumber
\end{eqnarray}
 
(\ref{lasteqn}) becomes
\begin{eqnarray}
\ & \ &\leq \sum_{k=1}^{\infty}\sum_{n=0}^{\infty}  (k+1)/2 \sqrt{t/ M}(k+2)/2 \sqrt{t/ M}
e^{-C_{0} (\frac{k}{2} \sqrt{t/M})^{\beta}} \nonumber \\
 \ &\ &\times \  e^{-\frac{(2n+1)^2}{(k+1)^2}M/2}
 \left[
\frac{1}{k^{4}} \left( \frac{\pi ^{3}(2n+1)^{3} t^{2}}{k^{2}} +
 3\pi (2n+1)t\right)+  \frac{(2n+1)\pi}{k^{2}}  \right.  \nonumber \\
\ &\ & + \left. \left( \frac{1}{k^{2}}(\frac{(2n+1)^{2}\pi ^{2}t
}{k^{2}}\  \ + \  1) \right) \right]\label{nksum}
\end{eqnarray}
We work on each sum in (\ref{nksum}) in $n$ separately, but they are all similar. We first consider the second term in
the square brackets in (\ref{nksum}). Now for $|x|<1$,
$$\sum_{n=0}^{\infty}n x^{n^2}\leq \sum_{n=0}^{\infty}n x^{n}=\frac{x}{(1-x)^{2}}.$$
Thus
\begin{eqnarray}
\sum_{n=0}^{\infty}  e^{-((2n+1)/(k+1))^{2} M/2 } ((2n+1)\pi )\leq \pi \frac{e^{-(1/(k+1))^{2} M/2 }}{(1-
e^{-(1/(k+1))^{2} M/2 }) ^{2}} \nonumber \\
\leq \pi \frac{1}{(1-e^{-(1/(k+1))^{2} M/2 }) ^{2}},\nonumber
\end{eqnarray}
Now since $C_{0}(\sqrt{t/M} 1/2)^{\beta} >1$ for $t$ large,
\begin{eqnarray}
 \ &\ &\sum_{k=1}^{\infty}    \frac{1}{k^{2}} \frac{(k+1)}{2} (k+2)/2( \sqrt{t/ M})^{2}\frac{e^{-C_{0} (\frac{k}{2}
\sqrt{t/M})^{\beta}}}{ (1-e^{-(1/(k+1))^{2} M/2 })^{2}}\nonumber \\
\ &\ & \leq ( t/ M)e^{-C_{0}( \frac{1}{2} \sqrt{t/M})^{\beta}} \nonumber \\
\ &\ & \times \ \sum_{k=1}^{\infty}
 \frac{(k+1)(k+2)}{4k^{2}} \frac{e^{-    (k^{\beta}-1)}}{ (1-e^{-(1/(k+1))^{2} M/2 }) ^{2}}.\nonumber
\end{eqnarray}
Now since
$$
(k+1)(k+2)\frac{e^{-    (k^{\beta}-1)}}{ (1-e^{-(1/(k+1))^{2} M/2 }) ^{2}}
$$
is bounded  uniformly for all $k$, actually this quantity tends to zero as $k$ tends to infinity, we get the desired conclusion. Alternately, for some $p>0$ and $C>0$ independent of $k$ we have
$$
(k+1)(k+2)\frac{e^{-    (k^{\beta}-1)}}{ (1-e^{-(1/(k+1))^{2} M/2 }) ^{2}}\leq k^{p}e^{-Ck^{\beta}}.
$$
Then since $\sum_{k=1}^{\infty}k^{p}e^{-Ck^{\beta}} < \infty$, we get the desired lower bound.
 All other terms are handled in a similar fashion; 
we use the following power series expansions for $|x|<1$
 $$\sum_{n=0}^{\infty}x^{n}=\frac{1}{(1-x)},$$
$$ \sum_{n=0}^{\infty}n^{2}x^{n}=\frac{x(1+x)}{(1-x)^{3}},$$ and $$\sum_{n=0}^{\infty}n^{3}x^{n}=\frac{x(1+4x+x^{2})}{(1-x)^{4}}.$$
\end{proof}

\begin{lemma}\label{lemmaA.3}
Let $\beta >0$. Let $\xi$ be a positive random variable which satisfies $-\log P[\xi > t] \sim c t^{\beta}$, as $t\rightarrow \infty $. Let $\xi _{1} $
and $\xi _{2}$ be independent copies of
$\xi$ and independent  of the Brownian motion $Y$. Then
\begin{eqnarray}
\ &\ & P[\eta _{(-\xi _{1}, \xi _{2})}>t]
  = P_{0}[\eta _{(-\xi _{1}, \xi _{2})}>t] \nonumber \\
 \ & \  &= \int_{0}^{\infty} \ \int_{0}^{\infty } \left( \frac{\partial }
{\partial u} \frac{\partial }{\partial v} P_{0} (\eta _{(-u,v)}>t) \right)  P[\xi
>u]P[\xi >v]dvdu .\nonumber
\end{eqnarray}
\end{lemma}
\begin{proof}
The proof follows from integration by parts in equation (\ref{eqint}), if we can show the
following two statements
\begin{equation}\label{endof1}
\lim_{u\rightarrow 0}P_{0}[\eta _{(-u,v)} > t] P[\xi > u]=\lim_{u\rightarrow \infty }P_{0}[\eta _{(-u,v)} > t] P[\xi > u]=0,
\end{equation}

\begin{equation}\label{endof2}
\lim_{v \rightarrow 0} (\frac{\partial }{\partial u} P_{0}[\eta _{(-u,v)} > t]) P[\xi > v]=\lim_{v\rightarrow \infty }
 (\frac{\partial }{\partial u}P_{0}[\eta _{(-u,v)} > t]) P[\xi > v]=0.
\end{equation}
The equation (\ref{endof1}) follows from the formula for  $P_{0}[\eta _{(-u,v)} > t]$  above  in equation (\ref{distr}) and the equation
(\ref{endof2}) is obvious from the asymptotics for $P[\xi > u]$ and equation (\ref{duu}).
\end{proof}

\textbf{Acknowledgments.} I would like to thank  Professor Rodrigo
 Ba\~{n}uelos, my academic advisor, for suggesting me this problem and for  his guidance on this paper. I also would like
 to thank Professor Dante DeBlassie and Professor Davar Khoshnevisan for useful comments on the early drafts of this work.


\begin{thebibliography}{99}

\bibitem{bds}  Ba\~{n}uelos, R., DeBlassie, R.D., Smits, R.: `The first
exit time of planar Brownian motion from the interior of a
parabola', \textit{Ann. Prob.} (2) \textbf{29} (2001), 882 - 901. 

\bibitem{vandenberg} van den Berg, M.: `Subexponential behavior of the Dirichlet heat kernel', \textit{J. Funct. Anal.} \textbf{198} (2003),
28 - 42. 

\bibitem{bgt} Bingham, N.H., Goldie, C.M., Teugels, J.L.: \textit{Regular Variation},
 Cambridge University Press, Cambridge, 1987.

\bibitem{burdzy1} Burdzy, K.: `Some path properties of iterated Brownian
motion', In Seminar on Stochastic Processes (E. \c{C}inlar, K.L.
Chung and M.J. Sharpe, eds.), \textit{Birkh\"{a}user}, Boston, (1993), 67 - 87.


\bibitem{burdzy2}  Burdzy, K.: `Variation of iterated
Brownian motion', In Workshops and Conference on Measure-valued
Processes, Stochastic Partial Differential Equations and
Interacting Particle Systems (D.A. Dawson, ed.)
 \textit{Amer. Math. Soc.} Providence, RI, (1994),35 - 53. 

\bibitem{bukh} Burdzy, K., Khoshnevisan, D.:
 `Brownian motion in a Brownian crack', \textit{Ann. Appl. Probabl.}
(3) \textbf{8} (1998), 708 - 748. 

\bibitem{csaki}  Cs\`{a}ki, E., Cs\"{o}rg\H{o}, M., F\"{o}ldes,  A., R\'{e}v\'{e}sz, P.: `The local time of iterated Brownian motion', \textit{J. Theoret. Probab.}
(3) \textbf{9} (1996), 717 - 743.

\bibitem{deblassie} DeBlassie, R. D.: `Iterated Brownian
motion in an open set', Preprint, 2004.

\bibitem{DSmits} DeBlassie, R.D., Smits, R.: `Brownian motion in twisted domains', Preprint, 2004.

\bibitem{eisenbaumshi} Eisenbum, N., Shi, Z.: `Uniform oscillations of the local time of iterated Brownian motion',
\textit{Bernoulli} (1) \textbf{5} (1999), 49 - 65. 

\bibitem{feller} Feller, W.: \textit{An Introduction to Probability
 Theory and its Applications}, Wiley, New York, 1971.

\bibitem{funaki}  Funaki, T.: `A probabilistic
 construction of the solution of some higher order parabolic differential
 equations',
\textit{Proc. Japan Acad. Ser. A. Math. Sci.} (5) \textbf{55} (1979), 176 - 179.


\bibitem{grad} Gradshteyn, I. S., Ryzhik, I. M.: \textit{Table of Integrals, Series and Products}, Academic, New
York, 1965.

\bibitem{kasahara}  Kasahara, Y.: `Tauberian theorems of exponential
type', \textit{J. Math. Kyoto Univ.} (2) \textbf{12} (1978), 209 - 219. 

\bibitem{klewis} Khoshnevisan, D., Lewis, T.M.: `Stochastic calculus for Brownian motion in a Brownian fracture',
\textit{Ann. Applied Probabl.} (3) \textbf{9} (1999), 629 - 667. 

\bibitem{lshi} Lifshits, M., Shi, Z.: `The first exit time of Brownian motion from a parabolic domain',
\textit{Bernoulli} (6) \textbf{8} (2002), 745 - 765. 

\bibitem{xiao}  Xiao, Y.: `Local times and related properties of multidimensional iterated Brownian
motion', \textit{J. Theoret. Probab.} (2) \textbf{11} (1998), 383 - 408. 

\end{thebibliography}
\end{document}